\documentclass[11pt,reqno]{amsart}
\usepackage{amssymb,verbatim,enumerate,ifthen}
\usepackage[mathscr]{eucal}
\usepackage[utf8]{inputenc}
\usepackage[T1]{fontenc}
\allowdisplaybreaks[1]
\def\N{\mathbb{N}}
\def\R{\mathbb{R}}

\newtheorem{theorem}{Theorem}
\newtheorem*{theorem*}{Theorem}
\def\Thm#1#2{\ifthenelse{\equal{#1}{*}}{\begin{theorem*}#2\end{theorem*}}
             {\begin{theorem}\label{T#1}#2\end{theorem}}}
\newtheorem{Atheorem}{Theorem}

\def\thm#1{Theorem~\ref{T#1}}

\newtheorem{proposition}[theorem]{Proposition}
\newtheorem*{proposition*}{Proposition}
\def\Prp#1#2{\ifthenelse{\equal{#1}{*}}{\begin{proposition*}#2\end{proposition*}}
             {\begin{proposition}\label{P#1}#2\end{proposition}}}

\newtheorem{corollary}[theorem]{Corollary}
\newtheorem*{corollary*}{Corollary}
\def\Cor#1#2{\ifthenelse{\equal{#1}{*}}{\begin{corollary*}#2\end{corollary*}}
             {\begin{corollary}\label{C#1}#2\end{corollary}}}

\newtheorem{conjecture}{Conjecture}
\newtheorem*{conjecture*}{Conjecture}
\def\Con#1#2{\ifthenelse{\equal{#1}{*}}{\begin{conjecture*}#2\end{conjecture*}}
             {\begin{conjecture}\label{Con#1}#2\end{conjecture}}}
\def\Cor#1{Conjecture~\ref{Con#1}}

\newtheorem{lemma}[theorem]{Lemma}
\newtheorem*{lemma*}{Lemma}
\def\Lem#1#2{\ifthenelse{\equal{#1}{*}}{\begin{lemma*}#2\end{lemma*}}
             {\begin{lemma}\label{L#1}#2\end{lemma}}}
\def\lem#1{Lemma~\ref{L#1}}
\newtheorem{Alemma}{Lemma}

\theoremstyle{definition}
\newtheorem{remark}[theorem]{Remark}
\newtheorem*{remark*}{Remark}
\def\Rem#1#2{\ifthenelse{\equal{#1}{*}}{\begin{remark*}\rm #2\end{remark*}}
             {\begin{remark}\label{R#1}\rm #2\end{remark}}}

\newtheorem{example}[theorem]{Example}
\newtheorem*{example*}{Example}
\def\Exa#1#2{\ifthenelse{\equal{#1}{*}}{\begin{example*}\rm #2\end{example*}}
             {\begin{example}\label{Ex#1}\rm #2\end{example}}}

\def\eq#1{{\rm(\ref{E#1})}}
\def\Eq#1#2{\ifthenelse{\equal{#1}{*}}
  {\begin{equation*}\begin{aligned}#2\end{aligned}\end{equation*}}
  {\begin{equation}\begin{aligned}\label{E#1}#2\end{aligned}\end{equation}}}

\def\sign{\mathop{\hbox{\rm sign}}\nolimits}

\def\comment#1{}

\begin{document}
\vspace{5mm}

\date{\today}

\title[]{Comparison of Gini means with fixed number of variables}

\author[R.\ Gr\"unwald]{Rich\'ard Gr\"unwald}
\author[Zs. P\'ales]{Zsolt P\'ales}
\address{Institute of Mathematics, University of Debrecen, H-4002 Debrecen, Pf. 400, Hungary}
\email{richard.grunwald@science.unideb.hu, pales@science.unideb.hu}

\thanks{The research of the first author was supported by the \'UNKP-23-3 New National Excellence Program of the Ministry for Innovation and Technology from the source of the National Research, Development and Innovation Fund. The research of the second author was supported by the K-134191 NKFIH Grant.}


\begin{abstract}	
	In this paper, we consider the global comparison problem of Gini means with fixed number of variables on a subinterval $I$ of $\R_+$, i.e., the following inequality
	\begin{align}\tag{$\star$}\label{ggcabs}
		G_{r,s}^{[n]}(x_1,\dots,x_n)
		\leq G_{p,q}^{[n]}(x_1,\dots,x_n),
	\end{align}
	where $n\in\N,n\geq2$ is fixed, $(p,q),(r,s)\in\R^2$ and $x_1,\dots,x_n\in I$.
	
	Given a nonempty subinterval $I$ of $\R_+$ and $n\in\N$, we introduce the relations
	\begin{align*}
		\Gamma_n(I)
		&:=\{((r,s),(p,q))\in\R^2\times\R^2\mid \eqref{ggcabs}\mbox{ holds for all }
		x_1,\dots,x_n\in I\}, \\
		\Gamma_\infty(I)
		&:=\bigcap_{n=1}^\infty\Gamma_n(I). 
	\end{align*}
	In the paper, we investigate the properties of these sets and their dependence on $n$ and on the interval $I$ and we establish a characterizations of these sets via a constrained minimum problem by using a variant of the Lagrange multiplier rule. We also formulate two open problems at the end of the paper.
\end{abstract}

\maketitle

\section{Introduction}

Throughout this paper, the symbols $\R$ and $\R_+$ will stand for the sets of real and positive real numbers, respectively, and $I$ will always denote a nonempty open real interval.

We begin by recalling the definition of the \emph{$n$-variable Gini mean corresponding to the parameters $(p,q)\in\R^2$}:
\Eq{*}{
	G_{p,q}^{[n]}(x_1,\dots,x_n)&:=\begin{cases}
		\bigg(\dfrac{x_1^p+\dots+x_n^p}{x_1^q+\dots+x_n^q}\bigg)^\frac{1}{p-q} & \mbox{if }p\neq q, \\[4mm]
		\exp\bigg(\dfrac{x_1^p\ln(x_1)+\dots+x_n^p\ln(x_n)}{x_1^p+\dots+x_n^p}\bigg) & \mbox{if } p=q,
	\end{cases}
	\qquad (x_1,\dots,x_n\in\R_+).
}
These means were invented by C.\ Gini in the paper \cite{Gin38}.
It is easy to observe that these means include the Hölder (or power) means. In particular, for all $p\in\R$, the Gini mean $G_{p,0}^{[n]}$ reduces to the $n$-variable $p$th power mean. According to a celebrated result of Aczél and Daróczy \cite{AczDar63c}, the homogeneous means among the so-called Bajraktarević means (\cite{Baj58,Baj63}) on the interval $\R_+$ are exactly the Gini means. 

The basic properties and identities for Gini means are summarized in the following assertion.

\Thm{BPG}{Let $(p,q)\in\R^2$ and $n\in\N$. Then
\begin{enumerate}
 \item $G_{p,q}^{[n]}=G_{q,p}^{[n]}$, i.e., Gini means are symmetric with respect to their parameters.
 \item $G_{p,q}^{[n]}\colon\R_+^n\to\R_+$ is a strict mean, i.e., for all $x_1,\dots,x_n\in\R_+$,
 \Eq{*}{
	\min\{x_1,\dots,x_n\}\leq G_{p,q}^{[n]}(x_1,\dots,x_n)\leq \max\{x_1,\dots,x_n\}
}
hold and both inequalities are strict if $\min\{x_1,\dots,x_n\}<\max\{x_1,\dots,x_n\}$.
 \item $G_{p,q}^{[n]}\colon\R_+^n\to\R_+$ is a homogeneous function, i.e., for all $t,x_1,\dots,x_n\in\R_+$,
 \Eq{*}{
   G_{p,q}^{[n]}(tx_1,\dots,tx_n)
   =tG_{p,q}^{[n]}(x_1,\dots,x_n).
 }
 \item $G_{p,q}^{[n]}\colon\R_+^n\to\R_+$ is an 
 infinitely many times differentiable symmetric function.
 \item For all nonzero $t\in\R$ and $x_1,\dots,x_n\in\R_+$,
 \Eq{*}{
 G_{tp,tq}^{[n]}(x_1,\dots,x_n)
 =\big(G_{p,q}^{[n]}(x_1^t,\dots,x_n^t)\big)^{\frac{1}{t}}.
 }
\end{enumerate}
}

Briefly, the aim of this paper is to investigate the global comparison problem of Gini means with fixed number $n$ of the variables from subinterval $I$ of $\R_+$, that is, to give necessary as well as sufficient conditions for the validity of the following inequality
\Eq{ggc}{
	G_{r,s}^{[n]}(x_1,\dots,x_n)
	\leq G_{p,q}^{[n]}(x_1,\dots,x_n),
}
where $n\in\N,n\geq2$ is fixed, $(p,q),(r,s)\in\R^2$ and $x_1,\dots,x_n\in I$.

\section{Preliminary results}

Given a nonempty subinterval $I$ of $\R_+$ and $n\in\N$, we introduce the relations
\Eq{*}{
  \Gamma_n(I)
  :=\{((r,s),(p,q))\in\R^2\times\R^2\mid \eq{ggc}\mbox{ holds for all }
  x_1,\dots,x_n\in I\} \qquad\mbox{and}\qquad
  \Gamma_\infty(I)
  :=\bigcap_{n=1}^\infty\Gamma_n(I). 
}
It is clear that $\Gamma_1(I)=\R^2\times\R^2$. The sets $\Gamma_2(I)$ and $\Gamma_\infty(I)$ have been characterized in the papers \cite{CziPal05,DarLos70,Pal88c,Pal89c,Pal92a}. To recall this 
result from those papers, we define the functions $\lambda,\mu\colon\R^2\to\R$ by
\Eq{*}{
  \lambda(u,v)
  :=\begin{cases}
     \min(u,v)& \mbox{if } u,v\geq0,\\[0.2mm] 
     0 & \mbox{if } uv<0,\\[0.2mm]
     \max(u,v)& \mbox{if } u,v\leq0, 
    \end{cases}
  \qquad\mbox{and}\qquad
  \mu(u,v)
  :=\begin{cases}
     \dfrac{|u|-|v|}{u-v} & \mbox{if } u\neq v,\\[2mm] 
     \sign(u) & \mbox{if } u=v.
    \end{cases}
}
For $(p,q)\in\R^2$, we also define the function $\chi_{p,q}\colon \R_+\to\R$ by
\Eq{*}{
	\chi_{p,q}(t)
	:=\begin{cases}
		\dfrac{t^p-t^q}{p-q} & \mbox{if } p\neq q, \\[3mm]
		t^p\ln(t) & \mbox{if } p=q.
	\end{cases}
}

\Thm{2,infty}{Let $I\subseteq\R_+$ be a subinterval and assume that $a:=\inf I<\sup I=:b$. 
If $a=0$ or $b=\infty$ hold, then
\Eq{*}{
  \Gamma_2(I)
  &=\big\{((r,s),(p,q))\in\R^2\times\R^2
  \mid r+s\leq p+q,\,\, 
  \lambda(r,s)\leq\lambda(p,q),\mbox{ \ and \ }
  \mu(r,s)\leq\mu(p,q)\big\},\\
 \Gamma_\infty(I)
  &=\big\{((r,s),(p,q))\in\R^2\times\R^2
  \mid \min(r,s)\leq\min(p,q)\mbox{ \ and \ }\max(r,s)\leq\max(p,q)\big\}.
}
If $0<a<b<\infty$ hold, then
\Eq{*}{
  \Gamma_2(I)
  &=\big\{((r,s),(p,q))\in\R^2\times\R^2
  \mid r+s\leq p+q
  \mbox{ \ and \ }
  G_{r,s}^{[2]}(a,b)\leq G_{p,q}^{[2]}(a,b)\big\},\\
 \Gamma_\infty(I)
  &=\big\{((r,s),(p,q))\in\R^2\times\R^2
  \mid 
  \chi_{r,s}\big(\tfrac{a}{b}\big)
  \leq\chi_{p,q}\big(\tfrac{a}{b}\big)
  \mbox{ \ and \ }
  \chi_{r,s}\big(\tfrac{b}{a}\big)
  \leq\chi_{p,q}\big(\tfrac{b}{a}\big)\big\}.
}}

The following lemma establishes the basic connection between the Gini mean $G_{p,q}^{[n]}$ and the corresponding function $\chi_{p,q}$.

\Lem{pq}{Let $(p,q)\in \R^2$, $n\in\N$ and let $\vartriangleleft$ denote any of the relations $<,\,\leq,\,=,\,\geq,\,>$. Then, for all $t,x_1,\dots,x_n\in\R_+$, the relation
\Eq{tG}{
  t \vartriangleleft G_{p,q}^{[n]}(x_1,\dots,x_n)
}
holds if and only if
\Eq{0chi}{
  0 \vartriangleleft \chi_{p,q}\Big(\frac{x_1}{t}\Big)
    +\dots+\chi_{p,q}\Big(\frac{x_n}{t}\Big)
}
is valid.}

\begin{proof}
Let us consider the case $p\neq q$ and assume first that $p>q$ holds. Let $\vartriangleleft$ denote any of the relations $<,\,\leq,\,=,\,\geq,\,>$. Then \eq{tG} has the form
\Eq{*}{
	t\vartriangleleft \bigg(\frac{x_1^p+\dots+x_n^p}{x_1^q+\dots+x_n^q}\bigg)^\frac{1}{p-q}\qquad(t,x_1,\dots,x_n\in\R_+).
}  
Taking the $(p-q)$th power, dividing by $t^p$ and multiplying by $(x_1^q+\dots+x_n^q)$ side by side, we arrive at
\Eq{*}{
	\bigg(\frac{x_1}{t}\bigg)^q+\dots+\bigg(\frac{x_n}{t}\bigg)^q\vartriangleleft\bigg(\frac{x_1}{t}\bigg)^p+\dots+\bigg(\frac{x_n}{t}\bigg)^p\qquad(t,x_1,\dots,x_n\in\R_+).
}
Subtracting $\big(\frac{x_1}{t}\big)^q+\dots+\big(\frac{x_n}{t}\big)^q$ from both sides and then dividing by $p-q$, which is a positive number in this case, we obtain \eq{0chi}. In the above calculation, all steps are easily reversible, so we have proved the equivalence of \eq{tG} and \eq{0chi} in this case. The proof in the case $p<q$ is similar, therefore it is omitted.

Finally, let us consider the case $p=q$. Then \eq{tG} has the form
\Eq{*}{
	t\vartriangleleft\exp\bigg(\frac{x_1^p\ln(x_1)+\dots+x_n^p\ln(x_n)}{x_1^p+\dots+x_n^p}\bigg)\qquad(t,x_1,\dots,x_n\in\R_+).
}
Taking the logarithm and then multiplying by $(x_1^p+\dots+x_n^p)$ side by side, we have
\Eq{*}{
	(x_1^p+\dots+x_n^p)\ln(t)\vartriangleleft x_1^p\ln(x_1)+\dots+x_n^p\ln(x_n)\qquad(t,x_1,\dots,x_n\in\R_+).
}
Subtracting $(x_1^p+\dots+x_n^p)\ln(t)$, applying a well-known identity for the logarithm function and then dividing by $t^p$, which is a positive number, we arrive at \eq{0chi}. In the above calculation, all steps are easily reversible, so we have proved the equivalence of \eq{tG} and \eq{0chi} in this case as well. \end{proof}

\section{Main results}

Our next result shows that $\big(\Gamma_n(I)\big)_{n\in\N}$ is a decreasing chain with respect to inclusion.

\Thm{Ch}{Let $I$ be a nonempty subinterval of $\R_+$ and $n,m\in\N$ with $n\leq m$. Then $\Gamma_m(I)\subseteq\Gamma_n(I)$.}

\begin{proof}
Let $n,m\in\N$ with $n\leq m$. For $n=m$, the statement is obvious, therefore we may assume that $n<m$. Let $((r,s),(p,q))\in\Gamma_m(I)$ be arbitrarily fixed. To show that $((r,s),(p,q))\in\Gamma_n(I)$, let $x_1,\dots,x_n\in I$ be arbitrary and, for $i\in\{n+1,\dots,m\}$, define
\Eq{*}{
  x_i:=t:=G_{r,s}^{[n]}(x_1,\dots,x_n).
}
According the \lem{pq}, the equality $t=G_{r,s}^{[n]}(x_1,\dots,x_n)$ yields that
\Eq{*}{
  0 = \chi_{r,s}\Big(\frac{x_1}{t}\Big)
    +\dots+\chi_{r,s}\Big(\frac{x_n}{t}\Big).
}
Since, for $i\in\{n+1,\dots,m\}$, we have that $\chi_{r,s}(x_i/t)=\chi_{r,s}(1)=0$, the above equality implies that
\Eq{*}{
  0 = \chi_{r,s}\Big(\frac{x_1}{t}\Big)
    +\dots+\chi_{r,s}\Big(\frac{x_m}{t}\Big).
}
Applying the \lem{pq} again, it follows that $t=G_{r,s}^{[n]}(x_1,\dots,x_m)$.

By the assumption $((r,s),(p,q))\in\Gamma_m(I)$, we have that
\Eq{*}{
  t=G_{r,s}^{[m]}(x_1,\dots,x_m)
  \leq G_{p,q}^{[m]}(x_1,\dots,x_m).
}
According to the \lem{pq}, the above inequality implies that
\Eq{*}{
  0 \leq \chi_{p,q}\Big(\frac{x_1}{t}\Big)
    +\dots+\chi_{p,q}\Big(\frac{x_m}{t}\Big).
}
Since, for $i\in\{n+1,\dots,m\}$, we have that $\chi_{p,q}(x_i/t)=\chi_{p,q}(1)=0$, the above inequality yields that
\Eq{*}{
  0 \leq \chi_{p,q}\Big(\frac{x_1}{t}\Big)
    +\dots+\chi_{p,q}\Big(\frac{x_n}{t}\Big).
}
In view of the \lem{pq}, it follows that $t\leq G_{p,q}^{[m]}(x_1,\dots,x_n)$, which obviously yields that 
\Eq{*}{
  G_{r,s}^{[n]}(x_1,\dots,x_n)
  \leq G_{p,q}^{[n]}(x_1,\dots,x_n).
}
is valid, which completes th proof of $((r,s),(p,q))\in\Gamma_n(I)$.
\end{proof}

To formulate the subsequent result, we introduce the following notations
\Eq{*}{
  tI:=\{tx\mid x\in I\} 
  \qquad
  I^\tau:=\{x^\tau\mid x\in I\},
  \qquad\mbox{and}\qquad
  I\cdot J:=\{xy\mid x\in I,\,y\in J\},
}
where $I,J\subseteq\R_+$ are nonempty subintervals, $t\in\R_+$ and $\tau\in\R\setminus\{0\}$. Clearly, $tI$, $I^\tau$ and $I\cdot J$ are also subintervals of $\R_+$. Furthermore, for $t\in\R_+$, one can easily establish the following equalities
\Eq{*}{
  \inf(tI)=t\inf(I), \qquad \inf(I^t)=\inf(I)^t,
   \qquad\mbox{and}\qquad
   \inf(I\cdot J)=\inf(I)\cdot\inf(J).
}

\Thm{bpg}{
	Let $I\subseteq\R_+$ be a nonempty subinterval and $n\geq2$. Then the following assertions hold.
\begin{enumerate}[(i)]
    \item If $J\subseteq\R_+$ is a nonempty subinterval and 
    $\inf(I\cdot I^{-1})\leq \inf(J\cdot J^{-1})$, then 
    \Eq{*}{
      \Gamma_n(I)\subseteq \Gamma_n(J).
    }
    In particular, if $I\supseteq J$ holds, then the above inclusion is valid. \\
    Furthermore, if $\inf(I\cdot I^{-1})=\inf(J\cdot J^{-1})$, then $\Gamma_n(I)=\Gamma_n(J)$. In particular, for all $t\in\R_+$, 
    \Eq{*}{
      \Gamma_n(tI)=\Gamma_n(I) \qquad\mbox{and}\qquad
      \Gamma_n(I^{-1})=\Gamma_n(I).
    }
    \item For all $t\in\R_+$, we have that 
    \Eq{*}{
       \Gamma_n(I^t)=t^{-1}\Gamma_n(I)
       \qquad\mbox{and}\qquad
	 \Gamma_n(I)
	 =\{((-p,-q),(-r,-s))\mid ((r,s),(p,q))\in\Gamma_n(I)\}.
	}	
	\item If $\inf(I\cdot I^{-1})=0$, then $\Gamma_n(I)$ is a cone (i.e., it is closed with respect to multiplication by positive numbers) and if $\inf(I\cdot I^{-1})>0$, then $\Gamma_n(I)$ is starshaped with respect to the point $((0,0),(0,0))$ (i.e., it is closed with respect to multiplication by numbers belonging to $[0,1]$). Furthermore,  
	\Eq{*}{
	  \Gamma_n(I)+(]-\infty,0]^2\times [0,\infty[^2)\subseteq\Gamma_n(I).
	}
	\item $\Gamma_n(I)$ is a nonempty closed subset of $\R^2\times\R^2$.
	\item $\Gamma_n(I)\circ\Gamma_m(I)\subseteq \Gamma_{\min(n,m)}(I)$ for every $n,m\in\N$.
\end{enumerate}
}

\begin{proof}
	(i) Assume that  $J\subseteq\R_+$ is a nonempty subinterval such that $\inf(I\cdot I^{-1})\leq \inf(J\cdot J^{-1})$. Let $((r,s),(p,q))\in\Gamma_n(I)$ be fixed. To show that $((r,s),(p,q))\in\Gamma_n(J)$ also holds, let first $x_1,\dots,x_n$ be arbitrary elements of the interior of $J$.
	Then, $x_1^{-1},\dots,x_n^{-1}$ are in the interior of $J^{-1}$. Consequently, for all $i\in\{1,\dots,n\}$, we have that $\inf(J)<x_i$ and $\inf(J^{-1})<x_i^{-1}$. Therefore, for all $i,j\in\{1,\dots,n\}$, the inequality $\inf(J\cdot J^{-1})=\inf(J)\cdot \inf(J^{-1})<x_i\cdot x_j^{-1}$ holds. According to our assumption, it follows that $\inf(I)\cdot \inf(I^{-1})=\inf(I\cdot I^{-1})<x_i\cdot x_j^{-1}$ is also valid for all $i,j\in\{1,\dots,n\}$. This implies that
	\Eq{*}{
	  \max_{1\leq i\leq n} \big(x_i^{-1}\cdot\inf(I)\big)
	  <\min_{1\leq j\leq n} \big(x_j^{-1}\cdot\sup(I)\big).
	}
	Choose $t>0$ such that 
	\Eq{*}{
	  \max_{1\leq i\leq n} \big(x_i^{-1}\cdot\inf(I)\big)
	  <t<\min_{1\leq j\leq n} \big(x_j^{-1}\cdot\sup(I)\big).
	}
	Then, these inequalities yield that, for all $i\in\{1,\dots,n\}$, we have $\inf(I)<tx_i<\sup(I)$, and hence $tx_i\in I$. Using that $((r,s),(p,q))\in\Gamma_n(I)$ holds,
	we can obtain that
	\Eq{*}{
	  G_{r,s}^{[n]}(tx_1,\dots,tx_n)
	  \leq G_{p,q}^{[n]}(tx_1,\dots,tx_n).
	}
	Using the homogeneity of Gini means, it follows that
	\eq{ggc} holds for arbitrary elements $x_1,\dots,x_n$ of the interior of $J$. Since the interior of $J$ is dense in $J$ and the Gini means are continuous functions, it follows that the above inequality is also valid for arbitrary elements $x_1,\dots,x_n$ of $J$. This proves that $((r,s),(p,q))\in\Gamma_n(J)$ and completes the proof of the first statement of assertion (i). 
	
	If $\inf(I\cdot I^{-1})=\inf(J\cdot J^{-1})$, then we have that $\inf(I\cdot I^{-1})\leq\inf(J\cdot J^{-1})$ and $\inf(J\cdot J^{-1})\leq\inf(I\cdot I^{-1})$ hold, whence according to the previous statement, we have $\Gamma_n(I)\subseteq \Gamma_n(J)$ and $\Gamma_n(J)\subseteq \Gamma_n(I)$, which show that the equality $\Gamma_n(I)\subseteq \Gamma_n(J)$ is valid, indeed.
	
	Let $t\in\R_+$ be fixed. Observe that with $J:=tI$, we have that $\inf(I\cdot I^{-1})=\inf(J\cdot J^{-1})$, therefore, $\Gamma_n(I)=\Gamma_n(J)$, i.e., $\Gamma_n(I)=\Gamma_n(tI)$ is true. Furthermore,
	observe that with $J:=I^{-1}$, we have that
	$\inf(J\cdot J^{-1})=\inf(I^{-1}\cdot (I^{-1})^{-1})=\inf(I^{-1}\cdot I)=\inf(I\cdot I^{-1})$.
	Therefore, $\Gamma_n(I)=\Gamma_n(J)$, i.e., $\Gamma_n(I)=\Gamma_n(I^{-1})$ is also true. Thus the proof of assertion (i) is complete.
	
	To verify the first statement of assertion (ii), let $t\in\R_+$ and $((r,s),(p,q))\in\Gamma_n(I^t)$. By definition, for all $x_1,\dots,x_n\in I^t$, we have that \eq{ggc} holds. Using the substitution $x_i:=u_i^t$, for all $u_1,\dots,u_n\in I$, it follows that
	\Eq{*}{
	  G_{r,s}^{[n]}(u_1^t,\dots,u_n^t)
	  \leq G_{p,q}^{[n]}(u_1^t,\dots,u_n^t).
	}
	This, in view of property (5) in Theorem 1, implies that, for all $u_1,\dots,u_n\in I$, 
	\Eq{*}{
	  G_{tr,ts}^{[n]}(u_1,\dots,u_n)
	  \leq G_{tp,tq}^{[n]}(u_1,\dots,u_n).
	}
	Therefore, $((tr,ts),(tp,tq))\in\Gamma_n(I)$, that is $((r,s),(p,q))\in t^{-1}\Gamma_n(I)$. Thus, we have proved the inclusion $\Gamma_n(I^t)\subseteq t^{-1}\Gamma_n(I)$.
	The proof of the reversed inclusion is analogous. 
	
	To show that the second statement of assertion (ii) is also valid, let $((r,s),(p,q))\in\Gamma_n(I)$. Then, for all $x_1,\dots,x_n\in I$, we have that \eq{ggc} holds. Using the substitution $x_i:=u_i^{-1}$, for all $u_1,\dots,u_n\in I^{-1}$, it follows that
	\Eq{*}{
	  G_{r,s}^{[n]}(u_1^{-1},\dots,u_n^{-1})
	  \leq G_{p,q}^{[n]}(u_1^{-1},\dots,u_n^{-1}).
	}
	This, in view of property (5) in Theorem 1, implies that, for all $u_1,\dots,u_n\in I^{-1}$, 
	\Eq{*}{
	  G_{-r,-s}^{[n]}(u_1,\dots,u_n)
	  \geq G_{-p,-q}^{[n]}(u_1,\dots,u_n).
	}
	Therefore, $((-p,-q),(-r,-s))\in\Gamma_n(I^{-1}) =\Gamma_n(I)$.
	
	To prove assertion (iii), observe that, for all $t\in\R_+$,
	\Eq{*}{
	  \inf(I^t\cdot (I^t)^{-1})
	  =\inf\big((I\cdot I^{-1})^t\big)
	  =\big(\inf(I\cdot I^{-1})\big)^t.
	}
	Therefore, if $\inf(I\cdot I^{-1})=0$, then $\inf(I^t\cdot (I^t)^{-1})=0$, which, according to assertion (i) implies that $\Gamma_n(I^t)=\Gamma_n(I)$.
	Now the first equality in assertion (ii) yields that 
	$\Gamma_n(I)=t\Gamma_n(I)$ for all $t\in\R_+$. Thus, $\Gamma_n(I)$ is a cone in this case.
	
	If $\inf(I\cdot I^{-1})>0$, then $\inf(I\cdot I^{-1})\in\,]0,1]$, consequently, for all $t\in\,]0,1]$,
	\Eq{*}{
	  \inf(I\cdot I^{-1})
	  \leq\big(\inf(I\cdot I^{-1})\big)^t
	  =\inf(I^t\cdot (I^t)^{-1}) .
	}
	Using assertions (i) and (ii), we can obtain that $\Gamma_n(I)\subseteq\Gamma_n(I^t)=t^{-1}\Gamma_n(I)$. Therefore, $t\Gamma_n(I)\subseteq\Gamma_n(I)$ for all $t\in\,]0,1]$, which proves that it is star-shaped with respect to the point $((0,0),(0,0))$. 
	
	To verify the last statement of assertion (iii), let $((r,s),(p,q))\in\Gamma_n(I)+\,]-\infty,0]^2\times [0,\infty[^2$ be arbitrary. Then there exists $((r',s'),(p',q'))\in\Gamma_n(I)$ such that $r\leq r'$, $s\leq s'$, $p'\leq p$, and $q'\leq q$. Then, according to \thm{2,infty}, it follows that, for all $x_1,\dots,x_n\in\R_+$,
	\Eq{*}{
	  G_{r,s}^{[n]}(x_1,\dots,x_n)
	  \leq G_{r',s'}^{[n]}(x_1,\dots,x_n) 
	  \qquad\mbox{and}\qquad
	  G_{p',q'}^{[n]}(x_1,\dots,x_n)
	  \leq G_{p,q}^{[n]}(x_1,\dots,x_n)
	}
	Using that $((r',s'),(p',q'))\in\Gamma_n(I)$, we also have that, for all $x_1,\dots,x_n\in I$,
	\Eq{*}{
	  G_{r',s'}^{[n]}(x_1,\dots,x_n)
	  \leq G_{p',q'}^{[n]}(x_1,\dots,x_n).
	}
	Combining these inequalities, we can conclude that
	for all $x_1,\dots,x_n\in I$,
	\Eq{*}{
	  G_{r,s}^{[n]}(x_1,\dots,x_n)
	  \leq G_{p,q}^{[n]}(x_1,\dots,x_n),
	}
	which shows that $((r,s),(p,q))\in\Gamma_n(I)$. Thus, we have completed the proof of assertion (iii).
	
	The nonemptiness of $\Gamma_n(I)$ follows from the inclusion $((p,q),(p,q))\in \Gamma_n(I)$. The closedness of $\Gamma_n(I)$ is an immediate consequence of the continuity of Gini means with respect to their parameters. Thus assertion (i) is shown.

	Finally, we prove assertion (v). Let $n,m\in\N$ and denote $k:=\min(n,m)$. Then, according to \thm{Ch}, 
	we have that $\Gamma_n(I)\cup \Gamma_m(I)\subseteq \Gamma_k(I)$. Therefore, $\Gamma_n(I)\circ\Gamma_m(I)\subseteq\Gamma_k(I)\circ\Gamma_k(I)$. Thus, it suffices to show that $\Gamma_k(I)\circ\Gamma_k(I)\subseteq\Gamma_k(I)$. 
	
	Let $((r,s),(p,q))\in\Gamma_k(I)\circ\Gamma_k(I)$. Then there exists $(u,v)\in\R^2$ such that $((r,s),(u,v))\in\Gamma_k(I)$ and $((u,v),(p,q))\in\Gamma_k(I)$. These inclusions imply that, 
	for all $x_1,\dots,x_k\in I$,
	\Eq{*}{
	 G_{r,s}^{[k]}(x_1,\dots,x_k)
	  \leq G_{u,v}^{[k]}(x_1,\dots,x_k) 
	  \qquad\mbox{and}\qquad
	  G_{u,v}^{[k]}(x_1,\dots,x_k)
	  \leq G_{p,q}^{[k]}(x_1,\dots,x_k).
    }
    Therefore, for all $x_1,\dots,x_k\in I$,
	\Eq{*}{
	 G_{r,s}^{[k]}(x_1,\dots,x_k)
	  \leq G_{p,q}^{[k]}(x_1,\dots,x_k),
    }
	which proves that $((r,s),(p,q))\in\Gamma_k(I)$.
\end{proof}

In the following result we characterize the elements of $\Gamma_n(I)$ via a conditional minimum problem.

\Thm{CEP}{Let $I\subseteq\R_+$ be a nonempty subinterval and $n\geq2$. Then $((r,s),(p,q))\in\Gamma_n(I)$ if and only if, for all $u_1,\dots,u_n\in\R_+$ with 
\Eq{CEP1}{
\chi_{r,s}(u_1)+\dots+\chi_{r,s}(u_n)=0
\qquad\mbox{and}\qquad 
\max(u_1,\dots,u_n)\inf (I\cdot I^{-1})\leq \min(u_1,\dots,u_n),
} 
the inequality
\Eq{CEP2}{
  0\leq \chi_{p,q}(u_1)+\dots+\chi_{p,q}(u_n)
}
holds.}

\begin{proof} Assume that $((r,s),(p,q))\in\Gamma_n(I)$ and let $u_1,\dots,u_n\in\R_+$ satisfy \eq{CEP1}. Then the inequality in \eq{CEP1} implies that
\Eq{*}{
  u_i\inf (I\cdot I^{-1})\leq u_j \qquad(i,j\in\{1,\dots,n\}).
}
Equivalently, this inequality can be rewritten as
\Eq{*}{
  u_j^{-1}\inf I\leq u_i^{-1}\sup I \qquad(i,j\in\{1,\dots,n\}),
}
whence the following inequality is obtained 
\Eq{*}{
  \max(u_1^{-1},\dots,u_n^{-1})\inf I\leq \min(u_1^{-1},\dots,u_n^{-1})\sup I.
}
Thus, there exists a value $t>0$ such that
\Eq{*}{
  \max(u_1^{-1},\dots,u_n^{-1})\inf I\leq t\leq \min(u_1^{-1},\dots,u_n^{-1})\sup I.
}
These inequalities imply that
\Eq{*}{
  \inf I\leq tu_i\leq \sup I \qquad(i\in\{1,\dots,n\}),
}
which yields that
\Eq{*}{
  tu_1,\dots,tu_n\in \overline{I}.
}
In view of the assumption $((r,s),(p,q))\in\Gamma_n(I)=\Gamma_n(\overline{I})$, we get that
\Eq{*}{
  G_{r,s}^{[n]}(tu_1,\dots,tu_n)\leq G_{p,q}^{[n]}(tu_1,\dots,tu_n).
}
Using the homogeneity if Gini means, we conclude that
\Eq{*}{
  G_{r,s}^{[n]}(u_1,\dots,u_n)\leq G_{p,q}^{[n]}(u_1,\dots,u_n).
}
Applying the equality in \eq{CEP1} and the \lem{pq}, we can see that $G_{r,s}^{[n]}(u_1,\dots,u_n)=1$. Therefore,
\Eq{*}{
  1\leq G_{p,q}^{[n]}(u_1,\dots,u_n),
}
which, by \lem{pq} again, implies that \eq{CEP2} is valid, indeed.

Conversely, let us first assume that, for all $u_1,\dots,u_n\in\R_+$ which satisfy \eq{CEP1}, the inequality \eq{CEP2} holds. To show that the inclusion $((r,s),(p,q))\in\Gamma_n(I)$ is valid, let $x_1,\dots,x_n\in I$ be arbitrary and denote 
\Eq{tGrs}{
	t:=G_{r,s}^{[n]}(x_1,\dots,x_n),
	\qquad u_1:=\frac{x_1}{t},\quad\dots,\quad u_n:=\frac{x_n}{t}.
}
Then, by the homogeneity of Gini means, we have that
\Eq{*}{
  1=G_{r,s}^{[n]}(u_1,\dots,u_n),
}
According to the \lem{pq}, it follows that $\chi_{r,s}(u_1)+\dots+\chi_{r,s}(u_n)=0$ holds and, for all $i,j\in\{1,\dots,n\}$, we have that $u_iu_j^{-1}=x_ix_j^{-1}\geq \inf(I\cdot I^{-1})$, which prove that $u_1,\dots,u_n$ satisfy the condition \eq{CEP1}.

Due to our assumption, we conclude that $u_1,\dots,u_n$ satisfy the inequality \eq{CEP2}. Applying \lem{pq}, it follows that
\Eq{*}{
	G_{r,s}^{[n]}(u_1,\dots,u_n)=1\leq G_{p,q}^{[n]}(u_1,\dots,u_n).
} 
By the definition of $u_1,\dots,u_n$ and by the homogeneity of Gini means, multiplying this inequality by $t$ side by side, we arrive at
\Eq{*}{
	G_{r,s}^{[n]}(x_1,\dots,x_n)\leq G_{p,q}^{[n]}(x_1,\dots,x_n),
}
which shows that $((r,s),(p,q))\in\Gamma_n(I)$.
\end{proof}

\Thm{M1}{Assume that $0<a:=\inf I<b:=\sup I<\infty$. Then $((r,s),(p,q))\in\Gamma_n(I)$ if and only if, for all $u_1,\dots,u_n\in\R_+$ such that \\ either
\Eq{MA1}{
  \begin{cases}
  G_{r,s}^{[n]}(u_1,\dots,u_n)=1,\\
  \max(u_1,\dots,u_n)a<\min(u_1,\dots,u_n)b, \\
  \mbox{$\exists\rho\in\R$ such that}\quad
  \chi'_{p,q}(u_k)+\rho\chi'_{r,s}(u_k)=0 \quad(k\in\{1,\dots,n\})
  \end{cases}
}
or
\Eq{MA2}{
  \begin{cases}
  G_{r,s}^{[n]}(u_1,\dots,u_n)=1,\\
  \max(u_1,\dots,u_n)a=\min(u_1,\dots,u_n)b,\\
  \mbox{$\exists\rho\in\R$ such that}\quad
  \chi'_{p,q}(u_k)+\rho\chi'_{r,s}(u_k)
  \begin{cases}
  \geq0 & \mbox{if }u_k=\min(u_1,\dots,u_n),\\
  =0 & \mbox{if } u_k\in\,]\!\min(u_1,\dots,u_n),\max(u_1,\dots,u_n)[\,,\\
  \leq0 & \mbox{if }u_k=\max(u_1,\dots,u_n)
  \end{cases}
  \end{cases}
}
or
\Eq{MA3}{
  \begin{cases}
  G_{r,s}^{[n]}(u_1,\dots,u_n)=1,\\
  \max(u_1,\dots,u_n)a=\min(u_1,\dots,u_n)b,\\
  \mbox{$\exists\rho\in\R\setminus\{0\}$ such that}\quad
  \rho\chi'_{r,s}(u_k)
  \begin{cases}
  \geq0 & \mbox{if }u_k=\min(u_1,\dots,u_n),\\
  =0 & \mbox{if } u_k\in\,]\!\min(u_1,\dots,u_n),\max(u_1,\dots,u_n)[\,,\\
  \leq0 & \mbox{if }u_k=\max(u_1,\dots,u_n)
  \end{cases}
  \end{cases}
}
the inequality 
\Eq{MB}{
  1\leq G_{p,q}^{[n]}(u_1,\dots,u_n)
}
holds.}

\begin{proof}  Assume first that $((r,s),(p,q))\in\Gamma_n(I)$
and let $u_1,\dots,u_n\in\R_+$ satisfy one of the conditions \eq{MA1} or \eq{MA2}. Then they also fulfill
\Eq{*}{
G_{r,s}^{[n]}(u_1,\dots,u_n)=1
\qquad\mbox{and}\qquad 
\max(u_1,\dots,u_n)a\leq \min(u_1,\dots,u_n)b.
}
Using \lem{pq} and that $\inf(I\cdot I^{-1})=a/b$, these conditions are equivalent to \eq{CEP1}. Therefore, according to \thm{CEP}, the inequality \eq{CEP2} holds, by \lem{pq} again, implies that \eq{MB} is also satisfied.

To prove the reversed implication, we consider the following constrained minimization problem:
\Eq{minpr}{
  \hbox{minimize} \quad \chi_{p,q}(u_1)+\dots+\chi_{p,q}(u_n)
  \qquad \hbox{subject to}\quad (u_1,\dots,u_n)\in U,
}
where
\Eq{*}{
	U:=\{(u_1,\dots,u_n)\in\R_+^n\mid \chi_{r,s}(u_1)+\dots+\chi_{r,s}(u_n)=0\quad\mbox{and}\quad au_j\leq bu_i \quad(i,j\in\{1,\dots,n\})\}.
}
We first show that the set of admissible points $U$ is compact.
It is clear that $U$ is closed, we need only to verify that $U$ is bounded. Let $(u_1,\dots,u_n)\in U$ be fixed. First observe that the inequalities $\min(u_1,\dots,u_n)\leq 1\leq \max(u_1,\dots,u_n)$ must be valid. Indeed, if $1<\min(u_1,\dots,u_n)$ were true, then, for all $i\in\{1,\dots,n\}$, we would have that $\chi_{r,s}(u_i)>0$, which contradicts the equality $\chi_{r,s}(u_1)+\dots+\chi_{r,s}(u_n)=0$. The inequality $\max(u_1,\dots,u_n)<1$ leads to a contradiction similarly. Using also the inequalities in the definition of $U$, it follows that, for all $i,j\in\{1,\dots,n\}$,
\Eq{*}{
  au_j\leq b\min(u_1,\dots,u_n)\leq b
  \qquad\mbox{and}\qquad
  a\leq\max(u_1,\dots,u_n)\leq bu_i,
}
which prove that $u_1,\dots,u_n\in \big[\frac{a}{b},\frac{b}{a}\big]$. Therefore, $U$ is bounded indeed.

Consequently, the minimization problem \eq{minpr} has a solution. 
To find these solutions, we use the Lagrange Multiplier Rule (for a smooth problem with one equality and $n^2$ inequality constraints). The Lagrange function of this problem is:
\Eq{*}{
  L\big(u_1,\dots,u_n,\lambda,\mu,(\nu_{i,j})_{i,j\in\{1,\dots,n\}}\big)
  =\lambda\sum_{i=1}^n \chi_{p,q}(u_i)+\mu\sum_{i=1}^n \chi_{r,s}(u_i)+\sum_{i=1}^n\sum_{j=1}^n\nu_{i,j}(au_j-bu_i).
}
According to the standard results, if $(u_1,\dots,u_n)$ is a solution to \eq{minpr}, then there exist real multipliers $\lambda,\mu,(\nu_{i,j})_{i,j\in\{1,\dots,n\}}$ (not all zero) such that
\Eq{*}{
  \lambda\geq0,\qquad \nu_{i,j}\geq0,\qquad \nu_{i,j}(au_j-bu_i)=0 \quad(i,j\in\{1,\dots,n\}),
}
and, for all $k\in\{1,\dots,n\}$, 
\Eq{LE}{
  \lambda\chi'_{p,q}(u_k)+\mu\chi'_{r,s}(u_k)
  +a\sum_{i=1}^n\nu_{i,k}-b\sum_{j=1}^n\nu_{k,j}=0.
}
There are two cases:

Case I: If $\max(u_1,\dots,u_n)a<\min(u_1,\dots,u_n)b$.
In this case $au_j<bu_i$ for all $i,j\in\{1,\dots,n\}$. Therefore, by the transversality condition, $\nu_{i,j}=0$ holds for all $i,j\in\{1,\dots,n\}$ and \eq{LE} simplifies to
\Eq{LE1}{
  \lambda\chi'_{p,q}(u_k)+\mu\chi'_{r,s}(u_k)=0 
  \qquad (k\in\{1,\dots,n\})
}
and $(\lambda,\mu)\neq(0,0)$. We now show that $\lambda\neq0$. 
If $\lambda=0$, then $\mu\neq0$, and the equality \eq{LE1} simplifies to
\Eq{*}{
  \chi'_{r,s}(u_k)=0 
  \qquad (k\in\{1,\dots,n\}).
}
Since the function $\chi'_{r,s}$ may possess at most one zero which must be different from $1$, we have that $u_1=\dots=u_n\neq 1$, which contradicts the inequalities $\min(u_1,\dots,u_n)\leq 1\leq \max(u_1,\dots,u_n)$. Therefore, $\lambda>0$ has to be valid and then the equality \eq{LE1} simplifies to
\Eq{*}{
  \chi'_{p,q}(u_k)+\frac{\mu}{\lambda}\chi'_{r,s}(u_k)=0
  \qquad (k\in\{1,\dots,n\}).
}
The equality in the definition of the set $U$, according to \lem{pq}, is equivalent to the condition $G_{r,s}(u_1,\dots,u_n)=1$.
Therefore, $u_1,\dots,u_n$ satisfy condition \eq{MA1} with $\rho=\mu/\lambda$.

Case II: If $\max(u_1,\dots,u_n)a=\min(u_1,\dots,u_n)b$.
In this case, denote 
\Eq{*}{
A:=\{i\in\{1,\dots,n\}\colon u_i=\min(u_1,\dots,u_n)\},\qquad
B:=\{j\in\{1,\dots,n\}\colon u_j=\max(u_1,\dots,u_n)\}.
}
Then, neither $A$ nor $B$ is empty and, for all $(i,j)\in A\times B$, the equality $au_j=bu_i$ holds. On the other hand, for all $(i,j)\not\in A\times B$, we have that $au_j<bu_i$, which, by the transversality condition implies that $\nu_{i,j}=0$ in these cases. Therefore, \eq{LE} simplifies to
\Eq{LE2}{
  \lambda\chi'_{p,q}(u_k)+\mu\chi'_{r,s}(u_k)
  +a\sum_{i\in A}\nu_{i,k}-b\sum_{j\in B}\nu_{k,j}=0 
  \qquad (k\in\{1,\dots,n\}).
}
If $(\lambda,\mu)=(0,0)$, then for $k\in A$, we get
\Eq{*}{
  \sum_{j\in B}\nu_{k,j}=0. 
}
Since the terms of this sum are nonnegative, it follows that $\nu_{k,j}=0$ for all $(k,j)\in A\times B$, which contradicts the nontriviality of the multipliers. 

If we set $k\in A$ in \eq{LE2}, then we get
\Eq{*}{
  \lambda\chi'_{p,q}(u_k)+\mu\chi'_{r,s}(u_k)
  =b\sum_{j\in B}\nu_{k,j}\geq0.
}
Arguing similarly for $k\in \{1,\dots,n\}\setminus(A\cup B)$ and $k\in B$, we can conclude that 
\Eq{LE3}{
  \lambda\chi'_{p,q}(u_k)+\mu\chi'_{r,s}(u_k)
  \begin{cases}
   \geq 0 & \mbox{if } k\in A,\\
   =0 & \mbox{if } k\in\{1,\dots,n\}\setminus(A\cup B),\\
   \leq 0 & \mbox{if } k\in B.
  \end{cases}
}
If $\lambda>0$, then we can see that $u_1,\dots,u_n$ satisfy condition \eq{MA2} with $\rho=\mu/\lambda$. In the case when 
$\lambda=0$, then $\mu\neq0$ and
\Eq{*}{
  \mu\chi'_{r,s}(u_k)
  \begin{cases}
   \geq 0 & \mbox{if } k\in A,\\
   \leq 0 & \mbox{if } k\in B,\\
   =0 & \mbox{if } k\in\{1,\dots,n\}\setminus(A\cup B).
  \end{cases}
}
Therefore, in this case, $u_1,\dots,u_n$ satisfy condition \eq{MA3} with $\rho=\mu$.

Consequently, by the assumption of the theorem, the inequality 
\eq{MB} is satisfied, equivalently, the inequality \eq{CEP2} holds.

By now, we have proved that for the solutions $(u_1,\dots,u_n)$ of the constrained minimization problem, the inequality \eq{CEP2} is valid. Therefore, the inequality \eq{CEP2} must be true for all $(u_1,\dots,u_n)\in U$. Consequently, in view of \thm{CEP}, we have that $((r,s),(p,q))\in\Gamma_n(I)$.
\end{proof}

\section{Conjectures}

We conclude this paper by formulating two conjectures that we have not been able to verify.

The first conjecture seems to be a natural extension of the formula for $\Gamma_2(I)$, if $0<\inf I<\sup I<\infty$.

\Con{1}{Let $n\geq3$ and let $I$ be a subinterval of $\R_+$ such that $0<a:=\inf I<\sup I=:b<\infty$. Then
\Eq{*}{
  \Gamma_n(I)
  &=\big\{((r,s),(p,q))\in\R^2\times\R^2
  \mid r+s\leq p+q,\,\,
  G_{r,s}^{[n]}(a,\dots,a,b)\leq G_{p,q}^{[n]}(a,\dots,a,b),\\ &\hspace{4.2cm}\mbox{ \ and \ }
  G_{r,s}^{[n]}(a,b,\dots,b)\leq G_{p,q}^{[n]}(a,b,\dots,b)\big\}.
}
In other words, $((r,s),(p,q))\in\Gamma_n(I)$ if and only if
\Eq{*}{
  r+s\leq p+q,\qquad
  G_{r,s}^{[n]}(a,\dots,a,b)\leq G_{p,q}^{[n]}(a,\dots,a,b),
  \qquad\mbox{and}\qquad
  G_{r,s}^{[n]}(a,b,\dots,b)\leq G_{p,q}^{[n]}(a,b,\dots,b).
}}

For the case when either $\inf I=0$ or the $\sup I=\infty$, we have the following conjecture, which is the limiting case of \thm{M1}.

\Con{M2}{Assume that $0=\inf I\cdot I^{-1}$ and let $n\geq3$. Then $((r,s),(p,q))\in\Gamma_n(I)$ if and only if, for all $u_1,\dots,u_n\in\R_+$ such that 
\Eq{*}{
  \begin{cases}
  G_{r,s}^{[n]}(u_1,\dots,u_n)=1,\\
  \mbox{$\exists\rho\in\R$ such that}\quad
  \chi'_{p,q}(u_k)+\rho\chi'_{r,s}(u_k)=0 \qquad(k\in\{1,\dots,n\})
  \end{cases}
}
the inequality 
\Eq{*}{
  1\leq G_{p,q}^{[n]}(u_1,\dots,u_n)
}
holds.}


\end{document}